\documentclass[11pt, a4paper, reqno]{amsart}

\usepackage{amscd}
\usepackage{dsfont}
\usepackage{amsmath, upgreek}
\usepackage{amstext}
\usepackage{amsthm}
\usepackage{amssymb}
\usepackage{amsopn}
\usepackage{amssymb}
\usepackage{amsxtra}
\usepackage{amsfonts}
\usepackage{fancyhdr}
\usepackage[mathcal]{euscript}
\usepackage[english]{babel}
\usepackage[latin1]{inputenc}
\usepackage{bbold}
\usepackage{setspace}
\usepackage[svgnames]{xcolor}
\usepackage{adjustbox}

\usepackage{tikz}
\usepackage{tikz-cd}
\usepackage{tkz-euclide}
\usepackage{euscript}
\usepackage{hyperref}
\hypersetup{
    colorlinks = true,
    linkcolor = RoyalBlue,
    citecolor = Green,    
}

\makeatletter
\g@addto@macro\normalsize{%
  \setlength\abovedisplayskip{10pt}
  \setlength\belowdisplayskip{10pt}
  \setlength\abovedisplayshortskip{5pt}
  \setlength\belowdisplayshortskip{8pt}
}

\newtheorem{thm}{Theorem}

\newtheorem{lem}[thm]{Lemma}

\theoremstyle{definition}

\newtheorem{rmk}[thm]{Remark}

\renewcommand{\epsilon}{\varepsilon}

\newcommand{\id}{\operatorname{id}\nolimits}

\newcommand{\Bigsum}[2]{\ensuremath{\mathop{\textstyle\sum}_{#1}^{#2}}}
\newcommand{\Bigprod}[2]{\ensuremath{\mathop{\textstyle\prod}_{#1}^{#2}}}

\newcommand{\ran}{\operatorname{ran}\nolimits}
\newcommand{\cok}{\operatorname{cok}\nolimits}

\newcommand{\ind}{\operatorname{ind}\nolimits}

\newcommand{\NN}{\mathbb{N}}

\newcommand{\CC}{\mathbb{C}}

\newcommand{\KK}{\mathbb{K}}
\newcommand{\EE}{\mathbb{E}}

\newcommand{\DD}{\mathbb{D}}

\newcommand{\diam}{\hfill\raisebox{-8pt}{\scalebox{0.9}{\rotatebox{45}{$\qed$}}}}

\newcommand{\HDLCS}{\mathsf{LCS}_{\mathsf{HD}}}

\newcommand{\FreDN}{\mathsf{Fre}_{\mathsf{(DN)}}}

\newcommand{\LB}{\mathsf{LB}}
\newcommand{\LF}{\mathsf{LF}}

\newcommand{\LBr}{\mathsf{LB}_{\mathsf{reg}}}

\newcommand{\LFr}{\mathsf{LF}_{\mathsf{reg}}}

\renewcommand{\sharp}{\texttt{\#}}

\newcommand{\hypref}[1]{\hyperref[R0]{\color{RoyalBlue}{\bf #1}}}

\usepackage{xstring}

	\usetikzlibrary{matrix,arrows,decorations,decorations.pathmorphing,positioning,decorations.pathreplacing,shapes}
	\tikzset{commutative diagrams/.cd, 
		mysymbol/.style = {start anchor=center, end anchor = center, draw = none}}

\DeclareMathSymbol{\sm}{\mathbin}{AMSa}{"39}
\DeclareMathSymbol{\shortminus}{\mathbin}{AMSa}{"39}

\DeclareMathSymbol{\shortminus}{\mathbin}{AMSa}{"39}

\usetikzlibrary{arrows.meta}

\usepackage{enumitem}

\newlist{compactitem}{itemize}{3} 
\setlist[compactitem,1]{leftmargin=26pt, noitemsep, topsep=0pt}

\newlist{myitemize}{itemize}{1}
\setlist[myitemize,1]{leftmargin=26pt, noitemsep, topsep=0pt}

\renewcommand{\tau}{\uptau}
\renewcommand{\sigma}{\upsigma}

\begin{document}
\allowdisplaybreaks
$ $
\vspace{-50pt}

\title{Exact sequences of regular LF-spaces}

\dedicatory{Dedicated to Dietmar Vogt on the occasion of his 85th birthday.}

\author{Jos\'e Bonet\hspace{0.75pt}$^{\text{1}}$ and Sven-Ake Wegner\hspace{0.5pt}$^{\text{2}}$}

\renewcommand{\thefootnote}{}
\hspace{-1000pt}\footnote{\hspace{5.5pt}2020 \emph{MSC Classification}: 46A13, 46M18, 18E05.\vspace{1pt}}


\hspace{-1000pt}\footnote{\hspace{5.5pt}\emph{Key words}: LF-space, inductive limit, preabelian category, (one-sided) exact category. \vspace{1pt}}

\hspace{-1000pt}\footnote{\hspace{0pt}$^{1}$\,Universitat Polit\`ecnica de Val\`encia, Instituto de Matem\'atica Pura y Aplicada, Camino de Vera, 46022 Valencia,\linebreak\phantom{x}\hspace{1.2pt}Spain, e-mail: jbonet@mat.upv.es.\vspace{1.6pt}}

\hspace{-1000pt}\footnote{\hspace{0pt}$^{2}$\,Corresponding author: University of Hamburg, Department of Mathematics, Bundesstra\ss{}e 55, 20146 Ham-\linebreak\phantom{x}\hspace{1.2pt}burg, Germany, phone: +49\,(0)\,40\:2395\hspace{1pt}-\hspace{1pt}25120, e-mail: sven.wegner@uni-hamburg.de.\vspace{1.6pt}}


\begin{abstract}
We show that the category of regular LF-spaces has all kernels but, unlike the category of regular LB-spaces, does not have all cokernels. Its topological exact structure, maximal exact structure and maximal deflation-exact structure form a strictly increasing chain of proper subclasses of the class of all kernel-cokernel pairs. To the best of the authors' knowledge, regular LF-spaces in particular provide the first natural example of a category that has all kernels but is not preabelian and whose maximal deflation-exact structure is neither two-sided exact nor equal to the class of all kernel-cokernel pairs.
\end{abstract}

\maketitle

\vspace{-10pt}

\section{Preliminaries}\label{SEC1}

This short note deals with (one-sided) exact categories formed by locally convex Hausdorff spaces, or lcHs for brevity. For unexplained notation we refer to \cite{BPC, Prosmans, MV, HKRW}. Recall that a lcHs $E$ is an LF-space if it can be written as the locally convex inductive (or co\hspace{1pt}-)limit of Fr\'echet spaces $E_1\hookrightarrow E_2\hookrightarrow\cdots$ with linear and continuous maps as structure maps which may be assumed to be injective without loss of generality. An LF-space $E=\ind_{n\in\NN}E_n$ is regular if every bounded subset of $E$ is contained and bounded in some $E_n$. It is well-known that an LF-space is regular if and only if it is locally (or Mackey) complete. We consider the categories
$$
\LFr\subset\LF\subset\HDLCS
$$
of regular resp.\ all LF-spaces as full subcategories of all lcHs equipped with linear and continuous maps as morphisms. It is well-known that $\HDLCS$ is quasiabelian and straightforward to check that $\LF$ and $\LFr$ are additive. Given a morphism $f\colon E\rightarrow F$ in $\LF$, its kernel is the inclusion $f^{-1}(0)^{\flat}=\ind_{n\in\NN}f^{-1}(0)\cap E_n\hookrightarrow E$ and its cokernel is the quotient map $F\rightarrow F/\overline{f(E)}$. In particular, $\LF$ is preabelian. In $\LFr$ the  kernel of a morphism $f\colon E\rightarrow F$ computes as in $\LF$. If the local completion $(F/\overline{f(E)})^{\sharp}$ is an LF-space, then $F\rightarrow(F/\overline{f(E)})^{\sharp}$ is the $\LFr$-cokernel of $f$. While this is always true if $E,F\in\LBr$, we show below that it does not hold in $\LFr$ and, even more, that in $\LFr$ there exist morphisms that do not possess a cokernel at all.

\smallskip

By \cite[Table on p.~2114]{DS12} the class of topologically exact sequences $\EE_{\mathsf{top},\,\LF}$ with all three terms being LF is an exact structure on $\LF$. Copying verbatim the proof of \cite[Prop 3.1.2]{DS16} one gets that the maximal exact structure $\EE_{\mathsf{max},\,\LF}$ consists of all short algebraically exact sequences
\begin{equation}\label{SES}
0\longrightarrow E\stackrel{f}{\longrightarrow}F\stackrel{g}{\longrightarrow}G\longrightarrow0
\end{equation}
in $\LF$ such that $f(E)\subseteq F$ is a well-located subspace; adapting their proof shows in particular that every cokernel in $\LF$ is semistable. Since the cokernels of $\LF$ are precisely the surjections, by \cite[Prop 9.2]{HR24} the maximal deflation-exact structure $\DD_{\mathsf{max},\,\LF}$ coincides with all kernel-cokernel pairs $\CC_{\mathsf{all},\,\LF}$ which are precisely the short algebraically exact sequences. From well-known examples based on D.~Vogt's seminal work \cite{Vogt92}, see also the proof of Theorem \ref{COR}(iv), it follows $\EE_{\mathsf{top},\,\LF}\subset\EE_{\mathsf{max},\,\LF}\subset\DD_{\mathsf{max},\,\LF}$ with proper inclusions.

\section{Results}

\begin{thm}\label{THM1} In the category $\LFr$ not every morphism has a cokernel. More precisely:\vspace{3pt}
\begin{compactitem}

\item[(i)] There is an LF-space $E=\ind_{n\in\NN}E_n$ such that $E^{\sharp}$ is not an LF-space, but equals $F^{\NN}$ for some complete LB-space $F$.\vspace{3pt}

\item[(ii)] In the situation of (i) write $E=Y/X$ with a complete LF-space $Y$ and a closed subspace $X\subseteq Y$. Then the inclusion $i\colon X^{\flat}\hookrightarrow Y$ has no cokernel in $\LFr$.
\end{compactitem}
\end{thm}
\begin{proof} (i) Let $F=\ind_{n\in\NN}F_n$ be a countable inductive limit of Banach spaces such that $F_n\subset F_{n+1}$ is a proper dense subspace for every $n\in\NN$, and such that $F$ is a complete lcHs; a concrete example is 
$$
\operatorname{s}'=\mathop{\ind}_{n\in\NN}F_n\,\text{ with } F_n=\bigl(\bigl\{\xi\in\KK^{\NN}\:\big|\:\|\xi\|_n^2=\Bigsum{j=1}{\infty}j^{-2n}|\xi_j|^2<\infty\bigr\},\|\cdot\|_n\bigr).
$$
Define $E_n:=F_n^{\NN}$ and $G:=F^{\NN}$ each with product topology and $E:=\ind_{n\in\NN}E_n$ . Then $E$ is by definition an LF-space, $G$ is complete and \cite[Prop 8.7.3]{BPC} yields further that $E\subseteq G$ is a topological subspace. By \cite[Prop 8.4.15]{BPC} the properties that we assumed for $F$ give us that $G=\Bigprod{\NN}{}F$ is not an LF-space. We claim that $E\subseteq G$ is locally dense: Given $x=(x_i)_{i\in\NN}\in G$ we have $x^{(k)}:=(x_1,\dots,x_k,0,\dots)\in E$ for every $k\in\NN$ and
$$
\lim_{k\rightarrow\infty}k(x-x^{(k)})=\lim_{k\rightarrow\infty}k(0,\dots,0,x_{k+1},x_{k+2},\dots)=0
$$
coordinate\hspace{1pt}-wise, i.e., in the topology of $G$. This means that every element of $G$ is a local limit point of $E\subseteq G$ and thus $G=E^{\sharp}$, see \cite[Chapter 5]{BPC}.

\medskip

(ii) For simplification write $X\equiv X^{\flat}$ and denote by $j\colon E\hookrightarrow E^{\sharp}$ the inclusion into the local completion. Assume $c\colon Y\rightarrow C$ is a cokernel of $i$ in $\LFr$. Then $ci=0$ and using that the quotient map $q\colon Y\rightarrow Y/X=E$ is a cokernel of $i$ in $\HDLCS$ we get a unique $f$ such that $fq=c$ holds. On the other hand, $c=\cok i$ holds in $\LFr$ and this yields for every $k\in\NN$ a map $g_k$ with $\pi_k j q=g_kc$, where $\pi_k\colon E^{\sharp}=\prod_{\NN}F\rightarrow F$ denotes the projection on the $k$-th factor. By the universal property of the product in $\HDLCS$ we get a unique $g$ with $\pi_k g=g_k$ for all $k\in\NN$ and thus the following diagram
\begin{equation*}
\begin{tikzcd}[column sep = 29pt, row sep = 14pt]
& & C \ar[rrdd,"g_k"]\ar[rdd,swap,"g"]\\
X\ar[r, hook, "i"]\ar[rru, bend left=20, "0"]\ar[rrd, bend right=20,swap, "0"] & Y\ar[ru,"c"]\ar[rd,swap, "q"] & \\
& & Y/X\ar[uu,"f"]\ar[r,hook, swap,"j"]&\prod_{\NN} F\ar[r,swap,"\pi_k"]&F
\end{tikzcd}
\end{equation*}
for which we claim $gf=j$. Indeed, we have $\pi_kjq=g_kc=\pi_kgfq$ and since $q$ is surjective $\pi_kj=\pi_kgf$. As this holds for every $k\in\NN$ we get $j=gf$. As $C$ is locally complete and $j$ is the inclusion of $E$ into its local completion we get a unique extension
$$
f^{\sharp}\colon E^{\sharp}=\Bigprod{\NN}{}F\longrightarrow C
$$
of $f$. Extending $j$ as well we obtain $gf^{\sharp}=j^{\sharp}=\id_{E^{\sharp}}$. Since $c$ is an $\LFr$-cokernel, it is an $\LFr$-epimorphism and thus has a dense range (use that $\KK\in\LFr$ and the Hahn-Banach Theorem). Consequently, $f(E)=f(q(Y))=c(Y)\subseteq C$ is dense and  therefore also $f^{\sharp}$ has a dense range. We may thus derive from $gf^{\sharp}=\id_{E^{\sharp}}$ that $f^{\sharp}g=\id_{C}$ holds as well. This would mean that $E^{\sharp}$ is an LF-space contradicting (i).
\end{proof}

The idea of the proof of Theorem \ref{THM1}(i) is in the same spirit as the proof of S.~Dierolf's \cite[Proposition in Section II]{Die81}.

\smallskip

To proceed, we need the following lemma whose LB-case is due to J.~Wengenroth, see \cite[Lem 4.12]{W25}.

\begin{lem}\label{LEM} Consider a short algebraically exact sequence of LF-spaces as in \eqref{SES} and assume that $f$ is a weak isomorphism onto its range. If $E$ and $G$ are regular then so is $F$.
\end{lem}
\begin{proof}Write $E=\ind_{n\in\NN}E_n$, $F=\ind_{n\in\NN}F_n$, $G=\ind_{n\in\NN}G_n$ and assume w.l.o.g.\ that $f$ is the inclusion of a well-located subspace induced by inclusions $f_n\colon E_n\rightarrow F_n$ and that $g$ is the quotient map induced by quotient maps $g_n\colon F_n\rightarrow G_n$.

\smallskip

Let $(x_k)_{k\in\NN}\subseteq F$ be a local Cauchy sequence. Then $(g(x_k))_{k\in\NN}\subseteq G$ converges locally to some $z\in G$, say in $G_B$ for a bounded disk $B\subseteq G$. By regularity there exists $n\in\NN$ such that $B$ is contained and bounded in $G_n$. Then $z$ also belongs to $G_n$ and $g(x_k)\rightarrow z$ converges in $G_n$. Using that $g_n\colon F_n\rightarrow G_n$ is open between metrizable spaces we find $y\in F_n$, $(y_k)_{k\in\NN}\subseteq F_n$ such that $g_n(y)=z$, $g_n(y_k)=g(x_k)$ and $y_k\rightarrow y$ in $F_n$. This implies $(x_k-y_k)_{k\in\NN}\subseteq g^{-1}(0)=f(E)$ and since $(y_k)_{k\in\NN}\subseteq F_n$ is in particular locally convergent, we get that $(x_k-y_k)_{k\in\NN}$ is locally Cauchy in $F$. We find a bounded disk $D\subseteq F$ such that $(x_k-y_k)_{k\in\NN}$ is contained and Cauchy in $F_D$. By intersecting with $f(E)$ we may assume that $D\subseteq (f(E),\uptau_F|_{f(E)})$ is bounded. Since $f$ is the inclusion of a linear subspace and a weak isomorphism onto its range we get that $D\subseteq E$ is bounded and a disk. Consequently, $(x_k-y_k)_{k\in\NN}$ is locally Cauchy in $E$, by regularity locally convergent in $E$ and therefore also locally convergent in $F$. Finally, $(x_k)_{k\in\NN}=((x_k-y_k)+y_k)_{k\in\NN}$ is locally convergent in $F$ as well. 
\end{proof}

Now we get the following result on (one-sided) exact structures on $\LFr$.

\begin{thm}\label{COR}The category $\LFr$ is additive and admits all kernels, while in it not all cokernels exist, i.e., it is not preabelian. Moreover the following holds: \vspace{2pt}
\begin{compactitem}

\item[(i)] $\EE_{\mathsf{top},\,\LFr}$ is an exact structure,\vspace{2pt}
\item[(ii)] $\EE_{\mathsf{max},\,\LFr}=\EE_{\mathsf{max},\,\LF}\cap\LFr$ is the maximal exact structure,\vspace{2pt}
\item[(iii)] $\DD_{\mathsf{max},\,\LFr}=\DD_{\mathsf{max},\,\LF}\cap\LFr$ is the maximal deflation-exact structure,\vspace{2pt}
\item[(iv)] $\EE_{\mathsf{top},\,\LFr}\subset\EE_{\mathsf{max},\,\LFr}\subset\DD_{\mathsf{max},\,\LFr}\subset\CC_{\mathsf{all},\,\LFr}$ are all proper subclasses.
\end{compactitem}
\end{thm}
\begin{proof}(i) Lemma \ref{LEM} implies in particular that $\LFr$ is extension closed in $(\LF,\EE_{\mathsf{top},\,\LF})$.

\smallskip

(iii) $\DD_{\mathsf{max},\,\LFr}$ consists of those kernel-cokernel pairs $E\stackrel{f}{\rightarrow}F\stackrel{g}{\rightarrow}G$ in which $g$ is a semistable cokernel, see \cite[Prop 9.2]{HR24}, and by observing that  \cite[Prop 2.2.3]{SiegThesis} holds also without the preabelianness assumed there, every semistable $\LFr$-cokernel is surjective. For the other inclusion observe that in $\LFr$ the pullback of a surjection along any morphism is again a surjection and hence a cokernel.

\smallskip

(ii) Lemma \ref{LEM} shows that $\LFr\subset(\LF,\EE_{\mathsf{max},\,\LF})$ is extension closed, hence $\EE_{\mathsf{max},\,\LF}\cap\LFr$ is an exact structure and thus contained in $\EE_{\mathsf{max},\,\LFr}$. For the other inclusion let $(f,g)$ as in \eqref{SES} be $\LFr$-stable. By the same argument as in (iii) we get that the latter is algebraically exact and it needs to be checked that $f(E)\subseteq F$ is well-located. By a lemma due to B.~Dierolf, see \cite[Lem 4.10]{W25}, this is equivalent to
$$
\forall\:t\colon E\rightarrow\KK\colon \KK\stackrel{q_{\KK}}{\longrightarrow}Q:=\textstyle\frac{F\hspace{1pt}\oplus\hspace{1pt}\KK}{\overline{\ran[f\;-t]^{\operatorname{T}}}}\text{ is injective,}
$$
where $q_{\KK}$ is the second component of the quotient map $q=[q_F\;q_{\KK}]$. To show this consider the $\LFr$-pushout of $f$ along $t$, which exists by assumption, and is the $\LFr$-cokernel of $\bigl[\begin{smallmatrix}\phantom{-}f\\-t\end{smallmatrix}\bigr]$. Let us denote this cokernel by $c=[c_F\;c_{\KK}]\colon F\oplus\KK\rightarrow C$. Using that $q=\cok\bigl[\begin{smallmatrix}\phantom{-}f\\-t\end{smallmatrix}\bigr]$ in $\HDLCS$ we get a unique map $i$ with $iq=c$ which means that the $\LFr$-pushout factors through the $\HDLCS$-pushout:
\begin{equation*}
\begin{tikzcd}
E\ar[r,"f"]\ar[d,swap,"t"]&F\ar[d,"q_F"]\\
\KK\ar[d,equals]\ar[r,"q_{\KK}"] & Q\ar[d,"i"] \\
\KK\ar[r,swap,"c_{\KK}"] & C. 
\end{tikzcd}
\end{equation*}
But now, since $f$ is $\LFr$-semistable, $c_{\KK}=iq_{\KK}$ is a kernel in $\LFr$ and hence injective. It follows that $q_{\KK}$ must also be injective.

\smallskip

(iv) Recycling of the (well-known) counterexamples mentioned in \cite[Ex 4.16 and Ex 4.14]{W25} provides that the first two inclusions are proper. Taking the standard resolution of any non-regular LB-space and appending the inclusion into the local completion yields an example witnessing that the last inclusion is proper.
\end{proof}

\begin{rmk} We conclude with the following remarks.\vspace{1pt}

\begin{compactitem}

\item[(i)] The category $(\LFr, \DD_{\mathsf{max},\,\LFr})$ is strongly deflation-exact and has kernels but not all kernels are inflations. It can thus be treated with the methods developed in \cite{HKRW,LW26} while results that assume two-sided exactness, preabelianness or admissible kernels do not apply.\vspace{2pt}

\item[(ii)] It is easy to see that there exist categories that have all kernels but are not preabelian and whose maximal deflation-exact structure is neither two-sided exact nor equal to the class of all kernel-cokernel pairs. Indeed, the product category $\FreDN\times\LBr$ provides an example, as follows by combining \cite[Rmk 4.8]{W25} and properties of $\FreDN$ which have been established jointly with J.~Wengenroth, see \cite[Section 7]{LW26}. The category $(\LFr, \DD_{\mathsf{max},\,\LFr})$ seems, to the best of the authors' knowledge, however to be the first natural example.\vspace{2pt}

\item[(iii)] The category of regular LB-spaces has cokernels since the local completion of an LB-space is always an LB-space. The (well-known) proof for the latter uses that LB-spaces are DF, see \cite[Lem 4.3]{W25}.\diam

\end{compactitem}
\end{rmk}

\section*{Acknowledgments}

\vspace{-2pt}

\small

The authors thank the referee for the prompt and careful review and for several helpful comments on the historically correct attribution of ideas and results used in this paper.

\normalsize

\vspace{-3pt}

\section*{Declaration of the Use of AI}

\vspace{-3pt}

\small

On 13 October 2026, the second author asked ChatGPT-6 Astra for an example in functional analysis in which every morphism has a kernel but which is not preabelian \emph{and} in which $\EE_{\mathsf{max}}\not=\DD_{\mathsf{max}}$ holds. It suggested $\FreDN\times\LB$, an example already known to him. The example presented here was conceived by the authors without AI assistance. ChatGPT subsequently assisted with manuscript revision, identified and corrected one mathematical error and shortened two arguments. The authors have carefully checked all its suggestions and take full responsibility for the results.

\normalsize

\vspace{-0pt}

\section*{Funding Information}

\vspace{-2pt}

\small

This research was part of project 507660524 funded by the German Research Foundation, project PID2024-162128NB-I00 funded by MICIU/AEI/10.13039/501100011033/FEDER, UE and project\linebreak{}CIAICO/2023/242 funded by Conselleria de Innovaci\'on, Universidades, Ciencia y Sociedad Digital de la Generalitat Valenciana.


\end{document}